\documentclass[preprint,11pt]{elsarticle}
\usepackage{amssymb,amsmath,amsthm}
\numberwithin{equation}{section}

\def\RR{\hbox{I\kern-.2em\hbox{R}}}

\usepackage{fullpage}
\newtheorem{guess}{Theorem}
\newtheorem{uess}{Lemma}
\newtheorem{corollary}{Corollary}
\newtheorem{definition}{Definition}
\newtheorem{example}{Example}
\newtheorem{remark}{Remark}

\begin{document}

\begin{frontmatter}

\title{On oscillation of difference equations with continuous time and variable 
delays}

\author[label1]{Elena Braverman}
\author[label1]{William T. Johnson}
\address[label1]{Dept. of Math. and Stats., University of
Calgary,2500 University Drive N.W., Calgary, AB, Canada T2N 1N4} 


\begin{abstract}
We consider existence of positive solutions for a difference equation with continuous time, variable coefficients and delays
$$
x(t+1)-x(t)+ \sum_{k=1}^m a_k(t)x(h_k(t))=0, 
\quad a_k(t) \geq 0, ~~h_k(t) \leq t, \quad t \geq 0, \quad k=1, \dots, m.
$$
We prove that for a fixed $h(t)\not\equiv t$, a positive solution may exist for $a_k$ exceeding any prescribed $M>0$, as well as for constant positive $a_k$ with $h_k(t) \leq t-n$, where $n \in {\mathbb N}$ is arbitrary and fixed. 
The point is that for equations with continuous time, non-existence of positive solutions with $\inf x(t)>0$ on any bounded interval should be considered rather than oscillation.   
Sufficient conditions when such solutions exist or do not exist are obtained.   We also present an analogue of the 
Gr\"{o}nwall-Bellman inequality for equations with continuous time, and 
examine the question when the equation has no positive non-increasing solutions. Counterexamples illustrate the role of variable delays.
\end{abstract}

\end{frontmatter}

\section{Introduction}
\label{sec1}

For both delay differential 
\begin{equation}
\label{01}
x^{\prime}(t) + \sum_{k=1}^m a_k(t) x(h_k(t)) = 0, \quad a_k(t) \geq 0, ~~k=1,\dots, m, ~ t \geq 0
\end{equation}
and difference equations
\begin{equation}
\label{02}
x(n+1)-x(n) + \sum_{k=1}^m a_k(n) x(h_k(n)) = 0, \quad a_k(n) \geq 0, ~~k=1,\dots, m,~n\in {\mathbb N}_0=0,1,\dots
\end{equation}
the following properties were used when studying oscillation.
\begin{enumerate}
\item
If $a_k \geq \alpha_k>0$, $k=1, \dots, m$, each non-oscillatory solution tends to zero:
$\displaystyle \lim_{t \to \infty} x(t) =0$.
\item
Any non-oscillatory solution is eventually monotone: 
positive solutions are eventually non-increasing, while negative solutions are non-decreasing.
\item
All positive solutions on any finite non-empty segment $[a,b]$ satisfy
$\displaystyle \min_{t \in [a,b]} x(t) >0$,
where in the case of difference equations the segment contains integers only.
\end{enumerate}

In the present paper, we illustrate that none of these properties, generally, is relevant for difference equations with continuous time.

The main object of the present paper is the equation with several variable delays 
\begin{equation}
\label{6}
x(t+1)-x(t)+ \sum_{k=1}^m a_k(t)x(h_k(t))=0, ~~
a_k(t) \geq 0, ~~h_k(t) \leq t, ~ t \geq 0, ~ k=1, \dots, m.
\end{equation}

\begin{definition}
A solution of \eqref{6} is {\bf oscillatory} if for any $t^{\ast}$ there exist $t_1>t^{\ast}$ and $t_2>t^{\ast}$
such that $x(t_1)<0<x(t_2)$. By a {\bf non-oscillatory} we mean an eventually positive or an eventually negative solution.
Equation~\eqref{6} is {\bf non-oscillatory} if it has a non-oscillatory solution and {\bf oscillatory} otherwise.
\end{definition}

For equations \eqref{01} and \eqref{02} with one delay term ($m=1$),
\begin{equation}
\label{04}
x^{\prime}(t) + a(t) x(h(t)) = 0, \quad a(t) \geq 0, \quad t \geq 0, 
\end{equation}
\begin{equation}
\label{05}
x(n+1)-x(n) + a(n) x(h(n)) = 0, \quad a(n) \geq 0, ~~~n=0,1,\dots
\end{equation}
the effect of non-monotonicity in variable delays leads to the conclusion that there is no $A>0$ and $B>0$ such that
$\displaystyle \limsup_{t \to \infty} \int_{h(t)}^t a(s)~ds>A$ implies oscillation of \eqref{04} or
$\displaystyle \limsup_{n \to \infty} \sum_{k=h(n)}^n a(k) >B$ guarantees oscillation of \eqref{05}, see \cite{ProcAMS,Karpuz2011}. 
Some mistakes in the previous results were reported in \cite{ProcAMS}.
However, if we replace $h$ with $g(t)=\sup_{s \leq t} h(s)$ and $g(n)=\sup_{k \leq n} h(k)$, respectively, $g$ becomes a non-decreasing function,
and any $A>1$ or $B>1$ in the above inequalities will imply oscillation \cite{GL}.

Some of oscillation and non-oscillation results for \eqref{01} and \eqref{02} were unified in the framework 
of equations on time scales, see, for example, \cite{Karpuz2010,Karpuz2008,Karpuz2016} and references therein.

Obviously, for \eqref{05} the inequality $\displaystyle \limsup_{n \to \infty} a(n)~>1$ leads to oscillation.
In this note, we prove that for a linear difference 
equation with continuous time and a variable argument
\begin{equation}
\label{1}
x(t+1)-x(t)+a(t)x(h(t))=0, \quad h(t) \leq t, \quad t \geq 0,
\end{equation}
no limitation $\displaystyle \limsup_{n \to \infty} a(n)~>M$ guarantees oscillation. Though in \eqref{1} the coefficient is also assumed to be variable, in the counterexamples we can even assume a constant $a(t)$. 

Sufficient oscillation conditions for the equation with constant delays and variable coefficients 
\begin{equation} \label{intro3}
x(t+1)-x(t)+\sum_{k=1}^n a_k(t)x(t-\sigma_k)=0, \ t \geq 0
\end{equation}
claim that, for high enough lower bound for $a_k$, all solutions oscillate \cite{Meng2005,Meng}. 
For the results on oscillation of \eqref{intro3}, its partial cases and some generalizations, 
see  \cite{Meng2005,Meng,Zhang} and references therein.

\begin{uess} \cite{Meng}
\label{lemMeng}
Let $a_k(t)$ be bounded, continuous, $\displaystyle \inf_{t\geq 0}a_k(t) \geq \underline{a_k}\geq 0$ and either
\begin{equation}
\label{eq1lemMeng}
\sum_{k=1}^n \underline{a_k} \frac{(\sigma_k+1)^{\sigma_k+1}}{\sigma_k^{\sigma_k}} > 1
\end{equation}
or
\begin{equation*}
n\left(\prod_{k=1}^n\underline{a_k}\right)^{\frac{1}{n}} \frac{(\sigma+1)^{\sigma+1}}{\sigma^{\sigma}} > 1, \mbox{\rm ~~where ~~} \sigma = \frac{1}{n}\sum_{k=1}^n \sigma_k.
\end{equation*}
Then all solutions of \eqref{intro3} oscillate.
\end{uess}

In particular,  by \eqref{eq1lemMeng} in Lemma~\ref{lemMeng}, 
$\displaystyle \inf_{t\geq 0} a(t)>1$ implies oscillation of  the equation
\begin{equation}
\label{intro1}
x(t+1)-x(t)+a(t)x(t-\sigma)=0, ~\sigma>0,~a(t)>0
\end{equation}
with a constant delay and a variable coefficient, for which the oscillation result can be stated as follows. 

\begin{uess} \cite{Indian}
\label{lemma1}
Let $a$ be continuous. 
If for $t$ large enough the inequality
\begin{equation}
\label{intro2}
a(t) \leq \frac{\sigma^{\sigma}}{(\sigma+1)^{\sigma+1}}
\end{equation}  
holds, \eqref{intro1} has a non-oscillatory solution. 
\end{uess}

Note that the result of Lemma~\ref{lemma1} 
was proven as early as in 1997, see \cite{Indian2}, 
under the additional assumption that $a(t)$ is Lipschitz continuous.


For \eqref{6} with constant delays $h_k(t)= t-\tau_k$, $\tau_k>0$, $k=1, \dots, m$,
the main tool to study oscillation was the Laplace Transform \cite{Meng2005,Meng}, and solutions were assumed to be piecewise continuous.
Here we formally do not imply any smoothness requirements on solutions but in all the examples we consider piecewise continuous solutions, 
following this tradition.

The paper is organized as follows. In Section~\ref{sec2} we illustrate first that no limitation on either magnitude of the coefficients or delays
guarantees oscillation of all solutions. However, if we introduce an additional restriction that on each finite segment, a solution has a positive lower bound, it is possible to obtain sufficient oscillation conditions.
Section~\ref{sec3} contains sufficient conditions when a non-oscillatory solution exists. In Section~\ref{sec4} we explore the problem of existence
of a positive non-increasing (or a negative non-decreasing) solution of \eqref{6}. Finally, Section~\ref{sec5} involves discussion and outlines some open problems.

\section{Oscillation}
\label{sec2}

First of all, let us note that $a(t) \geq a_0>0$ does not guarantee that all solutions of \eqref{1} tend to zero.

\begin{example}
\label{ex1}
Consider \eqref{1} 
with
$\displaystyle a(t)= \frac{1}{2+0.5^{\lfloor t \rfloor-2}} \geq \frac{1}{4}>0$, $\displaystyle h(t) = \lfloor t \rfloor - 0.5^{\lfloor t \rfloor} \left( 1- \{ t\} \right)$, ~$t \geq 1$,
where $\lfloor t \rfloor$ and $\{ t \}$ are the integer (the maximal integer not exceeding $t$) and the fractional parts ($\{ t\} = 
t- \lfloor t\rfloor$) of $t$, respectively.
Then $\displaystyle x(t)=\left(1+0.5^{\lfloor t \rfloor} \right)(1-\{t\})$ is a solution satisfying 
$x(n)=1+0.5^{\lfloor t \rfloor} >1$, thus $x(t)$ does not tend to zero.
To check that $x$ is a solution, we notice that $\lfloor h(t) \rfloor= \lfloor t \rfloor-1$, 
$\displaystyle 1-\{ h(t)\}= 0.5^{\lfloor t \rfloor} (1-\{ t\})$, 
$\displaystyle x(h(t))= \left( 1 + 0.5^{ \lfloor t \rfloor -1} \right) 0.5^{ \lfloor t \rfloor} (1-\{ t \})$,
and 
\begin{align*}
 & x(t+1)-x(t)+a(t)x(h(t)) \\  = & \left(1+0.5^{\lfloor t \rfloor +1} \right)(1-\{t\}) - 
\left(1+0.5^{\lfloor t \rfloor} \right)(1-\{t\})
+ \frac{1}{2+0.5^{ \lfloor t \rfloor -2}} \left(1+0.5^{ \lfloor t \rfloor -1} \right) 0.5^{ \lfloor t \rfloor} (1-\{t\}) \\  = &
(1-\{t\}) \left( 1+0.5^{ \lfloor t \rfloor +1} -1-0.5^{ \lfloor t \rfloor} +  0.5^{ \lfloor t \rfloor +1} 
\frac{1+0.5^{ \lfloor t \rfloor -1}}{1+0.5^{ \lfloor t \rfloor -1}} \right)
\\
 = & (1-\{t\})  \left( - 0.5^{ \lfloor t \rfloor +1} + 0.5^{ \lfloor t \rfloor +1} \right) =0.
\end{align*}
\end{example}

Next, let us note that if we consider piecewise continuous solutions, no conditions on the magnitude of $a(t)>0$
or on the delay $t-h(t)$ can exclude the possibility that a positive solution of \eqref{1} exists.

\begin{guess}
\label{th1}
For any $M>0$, there exist $a(t)>0$ and $h(t) \leq t$ such that 
\eqref{1} has a non-oscillatory solution, where $a(t) \geq M$
for any $t \geq 0$. For a fixed $a(t) \equiv a>0$ and any $n\in {\mathbb N}$, 
there exists $h(t) \leq t-n$ such that \eqref{1} has a non-oscillatory solution.
\end{guess}
{\bf Proof.} First, let $M>0$ be given. We can always find $k\in {\mathbb N}$ such that
$2^k \geq M$. Let $a(t) \equiv 2^k$. Denote
\begin{equation}
\label{2}
h(t)=\lfloor t \rfloor - 2^{-k-2} (1-\{ t\}) <t.
\end{equation}
Then $\lfloor h(t) \rfloor = \lfloor t \rfloor -1$, $\{ h(t)\} = 1-2^{-k-2}(1-\{t\})$.
It is easy to check that 
\begin{equation}
\label{3}
x(t)= 2^{- \lfloor t \rfloor } (1-\{ t\}) 
\end{equation}
is positive and satisfies \eqref{1} with $h(t)$ as in \eqref{2} and $a(t) \equiv 2^k$.

Next, let us fix $a(t)$, say $a(t) \equiv a= \frac{1}{4}$ and choose an arbitrary $n \in {\mathbb N}$. Denote
\begin{equation}
\label{4}
h(t)=\lfloor t \rfloor - n - 2^{-n} (1-\{ t\}) <t-n,
\end{equation}
then $\lfloor h(t) \rfloor = \lfloor t \rfloor -n-1$, $\{ h(t) \} =1-2^{-n}(1-\{ t \})$, and $x(t)$ defined in \eqref{3} is positive and satisfies \eqref{1} with $h(t)$ as in \eqref{4} and $a(t) \equiv \frac{1}{4}$,
which concludes the proof.
\qed

Thus, to study non-oscillation, we will only consider solutions that satisfy 
\begin{equation}
\label{5}
\inf_{t \in [c,d]} x(t) >0 \mbox{~~ for any finite~~} c<d.
\end{equation}

For equation \eqref{6}, we introduce lower and upper sequences
\begin{equation}
\label{7}
\underline{a_k}(n) = \inf_{s \in [n,n+1)} a_k(s), \quad \overline{h_k}(n)
= \sup_{s \in [n,n+1)} \left\lfloor h_k(s) \right\rfloor, \quad n\in {\mathbb N},~ k=1, \dots, m.
\end{equation}

\begin{guess}
\label{th2}
Let $a_k(t) \geq 0$, $k=1, \dots, m$. Assume that all solutions of the difference equation
\begin{equation}
\label{8}
y(n+1)-y(n) + \sum_{k=1}^m \underline{a_k}(n) y\left(\overline{h_k}(n) \right) =0
\end{equation}
oscillate. Then there is no positive solution of 
\eqref{6} satisfying \eqref{5}.
\end{guess}
{\bf Proof.} Assume the contrary that $x$ is a non-oscillatory solution of \eqref{6} satisfying \eqref{5}. 
First, let $x$ be positive. Denote the sequence
\begin{equation}
\label{9}
y(n)=\inf_{s \in [n,n+1)} x(s)>0.
\end{equation}
Let us note that, due to \eqref{5}, this is a positive decreasing sequence. If \eqref{5} fails, like in Example~\ref{ex1} and in 
the proof of Theorem~\ref{th1},
it can be a zero sequence.

In \eqref{7}, for $t \in [n,n+1)$ we have $\displaystyle a_k(t) \geq \underline{a_k}(n)$
and $\displaystyle h_k(t) < \overline{h_k}(n)$, leading to  
$$
x(h_k(t)) \geq y \left(\overline{h_k}(n) \right)>0, \quad t \in [n,n+1),
$$
and the solution $x$ satisfies
$$
x(t+1)=x(t)- \sum_{k=1}^m a_k(t)x(h_k(t)) \leq x(t) - \sum_{k=1}^m \underline{a_k}(n)
y \left(\overline{h_k}(n) \right), \quad t \in [n,n+1).
$$
By the definition of $y$ in \eqref{9}, this implies
$\displaystyle
y(n+1) \leq y(n)- \sum_{k=1}^m \underline{a_k}(n) y \left(\overline{h_k}(n) \right),
$
or
\begin{equation}
\label{10}
y(n+1)-y(n) + \sum_{k=1}^m \underline{a_k}(n) y \left(\overline{h_k}(n) \right) \leq 0.
\end{equation}
However, the existence of a positive solution of \eqref{10} is equivalent to the existence of a positive solution
of \eqref{8}, see, for example, \cite[Theorem 3.1]{ADSA2006}. 
The case of a negative solution is considered similarly, as the existence of a negative solution of the difference equation
$$
y(n+1)-y(n) + \sum_{k=1}^m \underline{a_k}(\alpha,n) y \left(\overline{h_k}(n) \right) \geq 0
$$
implies the existence of a non-oscillatory solution to \eqref{8} \cite[Theorem 3.1]{ADSA2006}.
The contradiction shows that all solutions of \eqref{6} satisfying \eqref{5} oscillate.
\qed

\begin{corollary}
\label{cor_compare}
Let $a_k(t) \geq 0$, $k=1, \dots, m$. Assume that 
\begin{equation*}
a_k(t) \geq b_k(n)>0,~~ h_k(t) \leq g_k(n) \leq n,~~t\in [n,n+1),~~k=1, \dots, m,
\end{equation*}
and all solutions of the difference equation
\begin{equation*}
y(n+1)-y(n) + \sum_{k=1}^m b_k(n) y\left(g_k(n) \right) =0
\end{equation*}
oscillate. Then there is no positive solution of
\eqref{6} satisfying \eqref{5}.
\end{corollary}

Making a shift from $n\in {\mathbb N}$ to $n+\alpha$ for some $\alpha \in [0,1)$, we generalize Theorem~\ref{th2}.

\begin{guess}
\label{th2a}
Let $a_k(t) \geq 0$, $k=1, \dots, m$. Assume that for some $\alpha \in [0,1)$, all solutions of the difference equation
\begin{equation*}
y(n+1)-y(n) + \sum_{k=1}^m \underline{a_k}(\alpha,n) y\left(\overline{h_k}(\alpha,n) \right) =0
\end{equation*}
oscillate, where
$\displaystyle \underline{a_k}(\alpha, n)= \inf_{s \in [\alpha+n,\alpha+n+1)} a_k(s),~
\overline{h_k}(\alpha, n) = \sup_{s \in [\alpha+ n,\alpha+ n+1)} \left\lfloor h_k(s) \right\rfloor,
~ n\in {\mathbb N},~ k=1, \dots, m.$

Then there is no positive solution of \eqref{6} satisfying \eqref{5}.
\end{guess}

\begin{example}
Consider the equation
\begin{equation}
\label{ex1eq1}
x(t+1)-x(t) + \sum_{k=1}^{10} (0.03+0.003 \cos(kt)) x(t-2.1+0.1 \sin(kt))=0.
\end{equation}
Here $\underline{a_k}$ and $\overline{h_k}(n)$ introduced in \eqref{7} are ~
$\displaystyle
\underline{a_k}(n) \geq 0.027, ~ \overline{h_k}(n) \leq n-1, ~ n\in {\mathbb N},~ k=1, \dots, 10$
and ~ $\displaystyle \sum_{k=1}^{10} \underline{a_k}(n) \geq 0.27$.
Since $\displaystyle 0.27 > \frac{1}{4} =\frac{1^1}{2^2}$, the equation
$$
x(n+1)-x(n)+0.27x(n-1)=0
$$
is oscillatory \cite[Theorem 7.2.1]{GL}. By Corollary~\ref{cor_compare}, there are no positive solutions of \eqref{ex1eq1} satisfying \eqref{5}.
\end{example}

\section{Existence of positive solutions}
\label{sec3}

Further, consider existence of non-oscillatory (positive) solutions.
To this end, extend the definitions in \eqref{7} to lower and upper sequences for both coefficients and arguments
\begin{equation}
\label{13}
\underline{a_k}(n) = \inf_{s \in [n,n+1)} a_k(s), \quad 
\overline{a_k}(n) = \sup_{s \in [n,n+1)} a_k(s), \quad n\in {\mathbb N},~ k=1, \dots, m
\end{equation}
\begin{equation}
\label{14}
\underline{h_k}(n)
= \inf_{s \in [n,n+1)} \left\lfloor h_k(s) \right\rfloor, \quad
\overline{h_k}(n)
= \sup_{s \in [n,n+1)} \left\lfloor h_k(s) \right\rfloor, \quad n\in {\mathbb N},~ k=1, \dots, m.
\end{equation}

First of all, let us notice that existence of a positive solution of \eqref{intro1},
even with a constant coefficient
\begin{equation}
\label{14a}
x(t+1)-x(t)+q x(t-\sigma)=0, ~\sigma>0
\end{equation}
does not imply that the solution of the equation 
\begin{equation}
\label{14ab}
x(t+1)-x(t)+q x(h(t))=0, ~ t- \sigma \leq h(t) \leq t
\end{equation}
with the same initial conditions is positive.

\begin{example}
The equation 
$$
x(t+1)-x(t)+\frac{1}{4} x(t-1)=0
$$
has a positive solution $x(t)=2^{-t}$. Consider 
$$
x(t+1)-x(t)+\frac{1}{4} x(h(t))=0, ~~ h(t)= \left\{ \begin{array}{cc} t, & t \in [n,n+\frac{1}{2}),~n \in {\mathbb N}, \\
t-1,  & t \in [n+\frac{1}{2},n+1),~1 \leq n \leq 5,~ n \in {\mathbb N}, \\
t - \frac{1}{2},  & t \in [n+\frac{1}{2},n+1),~n \geq 6, ~n \in {\mathbb N}, \end{array}   \right.
$$
with the initial condition $x(t)=2^{-t}$, $t \in [-1,1)$.
The solution on $[0,6.5)$ is
$$
x(t)= \left\{ \begin{array}{cc} (\frac{3}{2})^n 2^{-t}, & t \in [n,n+\frac{1}{2}),~t<6.5, \vspace{2mm} \\
2^{-t},  & t \in [n+\frac{1}{2},n+1),~ t<6, \end{array}  \right.
$$
while on $t \in [6.5,7)$ we have
$\displaystyle x(t)= 2^{-t+1} - \frac{1}{4} \left( \frac{3}{2} \right)^5 2^{-t+\frac{3}{2}}<0, ~~ t 
\in [6.5,7).
$
\end{example}

Further, we present sufficient conditions for existence of a positive solution  
of \eqref{6} satisfying \eqref{5}. Similar conditions can be obtained for a negative solution with a negative supremum on any finite segment.

\begin{guess}
\label{th3}
Suppose that there exist positive non-increasing sequences 
$u(n)$ and $V(n)$ for $n \in {\mathbb Z}$ 
such that
\begin{equation}
\label{15}
u(n) \leq V(n), \quad n \geq 0,
\end{equation}
\begin{equation}
\label{16}
u(n+1) \leq u(n) - \sum_{k=1}^m \overline{a_k}(n) V\left( \underline{h_k}(n) \right), \quad n \geq 0,
\end{equation}
\begin{equation}
\label{17}
V(n+1) \geq V(n) - \sum_{k=1}^m \underline{a_k}(n) u\left( \overline{h_k}(n) \right), \quad n \geq 0.
\end{equation}  
Then there exists a positive solution $x(t)$ of
\eqref{6} satisfying \eqref{5}. Moreover, there is $x(t)$ for which
\begin{equation}
\label{18}
u(n)\leq x(t) \leq V(n), \quad t \in [n,n+1).
\end{equation}
\end{guess}

{\bf Proof.} Let us consider the solution of \eqref{6} with the initial function
$x(t)=\varphi(\lfloor t \rfloor)$, $t<1$, where $u(n) \leq \varphi(n) \leq V(n)$ for $n \leq 0$. 
In fact we are going to verify $u(n)\leq x(t) \leq V(n)$, $t \in [n,n+1)$,  which
would imply \eqref{18} and the statement of the theorem.

We prove \eqref{18} by induction. We have $\varphi(n)=u(n)=x(t)=V(n)$ for any $t\in [n,n+1)$, where $n \leq 0$.  
Further, for $t\in [0,1)$ and $t+1\in [1,2)$, we have by \eqref{16}
\begin{eqnarray*}
x(t+1) & = & x(t) - \sum_{k=1}^m a_k(t)x(h_k(t))= \varphi(0)-\sum_{k=1}^m a_k(t) \varphi(\lfloor h_k(t) \rfloor)
\\ & \geq & u(0)-\sum_{k=1}^m a_k(t) V(\lfloor h_k(t) \rfloor \geq 
u(0) - \sum_{k=1}^m \overline{a_k}(t) V( \underline{h_k}(0)) \geq u(1),
\end{eqnarray*}
and by \eqref{17}
\begin{eqnarray*}
x(t+1) & = & x(t) - \sum_{k=1}^m a_k(t)x(h_k(t)) = \varphi(0)-\sum_{k=1}^m a_k(t)\varphi(\lfloor h_k(t) \rfloor)
\\ & \leq & V(0) -\sum_{k=1}^m a_k(t) u(\lfloor h_k(t) \rfloor) \leq  V(0) -\sum_{k=1}^m \underline{a_k}(t) u(\overline{h_k}(0)) \leq 
V(1).
\end{eqnarray*}
Thus 
$u(0) \leq x(t) \leq V(0)$ for $t \in [1,2)$.

Further, let us assume that
$u(n) \leq x(t) \leq V(n)$ for $t \leq n+1$.
Then, for $t\in [n,n+1)$, $t+1 \in [n+1,n+2)$, by \eqref{17} we have
$$x(t+1) = x(t) - \sum_{k=1}^m a_k(t)x(h_k(t)) \leq V(n) - \sum_{k=1}^m \underline{a_k}(n) u\left( \overline{h_k}(n) \right) \leq V(n+1),$$
while \eqref{16} implies
$$x(t+1) = x(t) - \sum_{k=1}^m a_k(t)x(h_k(t)) \geq u(n) - \sum_{k=1}^m \overline{a_k}(n) V\left( \underline{h_k}(n) \right) \geq u(n+1),$$
therefore $u(n+1) \leq x(t) \leq V(n+1)$, $t \in [n+1,n+2)$, which concludes the induction step and the proof.
\qed

Consider \eqref{6} with piecewise constant coefficients and arguments.

\begin{corollary}
\label{cor1}
Suppose that 
$a_k(t) \equiv a_k(n)$,  $t \in [n,n+1)$, $\underline{h_k}(n)=\overline{h_k}(n)=\theta_k(n)$ for 
any $n \in {\mathbb N}$,
and the difference equation
\begin{equation}
\label{22}
z(n+1)-z(n)+\sum_{k=1}^m a_k(n) z(\theta_k(n))=0
\end{equation}
has a positive solution.
Then there exists a positive solution $x(t)$ of
\eqref{6} satisfying \eqref{5}. 
\end{corollary}
{\bf Proof.} Let $z(n)$ be a positive solution of \eqref{22}.
The sequence $z(n)$ satisfies \eqref{15} with $\varphi(n)=z(n)$ for $n \leq 0$,  and $u(n)=z(n)$, $V(n)=z(n)$, $n \geq 0$ satisfy \eqref{16},\eqref{17} with the equality signs. The application of Theorem~\ref{th3} 
concludes the proof.
\qed

Let us also note that under the conditions of the corollary, the positive solution of \eqref{6}
discussed in the proof is piecewise constant. Also, in Theorem~\ref{th3},  the initial point $t=1$ can be substituted with any~$t_0$.

%
%

\begin{example}
For the equation
\begin{equation}
\label{ex3eq1}
x(t+1)-x(t)+  0.5^{\lfloor t \rfloor+2} x\big(\lfloor t \rfloor-1-0.8 \cos t\big)=0
\end{equation}
$\underline{h}(n)=n-2$, $\overline{h}(n)=n$, $\overline{a}(n) = \underline{a}(n)= 0.5^{n+2}$. Denote
$V(n)=1$, $u(n)=0.5+0.5^{n+1}$, $n \geq 0$ and $u(n)=V(n)=1$, $n<0$.
Then \eqref{17} is obviously satisfied, while \eqref{16}  has the form
$$u(n)-u(n+1) = 0.5+0.5^{n+1}-\left(0.5+0.5^{n+2} \right) =  0.5^{n+2} \geq 0.5^{n+2} \cdot 1 = \overline{a}(n) V(\underline{h}(n)),$$
which is also true. Thus \eqref{ex3eq1} has a positive solution satisfying \eqref{5}.
\end{example}


\section{Existence of positive non-increasing solutions}
\label{sec4}

Finally, consider the problem of existence of positive non-increasing and negative non-decreasing solutions.
Here we consider variable and, generally, non-monotone arguments $h_k(t)$, however, we introduce non-decreasing functions
\begin{equation}
\label{phi_k}
g_k(t) = \sup_{s \leq t} h_k(s), \quad k=1, \dots , m,
\quad g(t) = \max_{1 \leq k \leq m} g_k(t).
\end{equation}
For any two real numbers $s \leq t$ we introduce two functions, each describing a finite set of numbers
\begin{equation}
\label{sets}
\displaystyle {\cal N}(s,t) = \left\{ \left. s+j \right| j=0,1,2, ~ s+j \leq t-1  \right\}, 
~
\displaystyle {\cal M}(s,t) = \left\{ \left. t-j \right| j \in {\mathbb N}, ~ t-j \geq s \right\}.
\end{equation}
Obviously, if $t=s+n$, $n\in {\mathbb N}$, these sets coincide ${\cal N}(s,t) = {\cal M}(s,t) = \{s, s+1, \dots s+n-1\}$.
For example, if $t=s+n+0.6$, $n\in {\mathbb N}$, we have ${\cal N}(s,t)= \{s, s+1, \dots s+n-1\}$, while
${\cal M}(s,t)=\{s+0.6, s+1.6, \dots s+n-0.4 \}$.
Further, we assume that the sum with no terms equals zero, while the product with no factors is one.
The following lemma evaluates the rate of the minimal decay of positive solutions. It can be treated as an analogue of the Gr\"{o}nwall-Bellman inequality for equations with continuous time.

\begin{uess}
\label{lem_flow}
Let $x(t)$ be a positive non-increasing solution of \eqref{6}. Then, for any $0\leq s < t < \infty$,
\begin{equation}
\label{G1}
x(t) \leq x(s) \prod_{r_j \in {\cal N}(s,t)} \left( 1- \sum_{k=1}^m a_k(r_j) \right),
~~ x(t) \leq x(s) \prod_{r_j \in {\cal M}(s,t)} \left( 1- \sum_{k=1}^m a_k(r_j) \right).
\end{equation}
\end{uess}
{\bf Proof.}
First, let us note that for $0\leq s <t<s+1$, under the assumption that $x$ is non-increasing, 
the inequalities in \eqref{G1} take an obvious form $x(t) \leq x(s)$.
Further, since $h_k(s)\leq s$ and $x$ is non-increasing, we have
$$
x(s+1)=x(s) - \sum_{k=1}^m a_k(s) x(h_k(s)) \leq x(s) - \sum_{k=1}^m a_k(s) x(s) = x(s) \left( 1- \sum_{k=1}^m a_k(s) \right).
$$
For any $t\in [s+1,s+2)$, the set ${\cal N}$ consists of $r_1=s$ only, therefore
$$x(t) \leq x(s+1) \leq x(s) \prod_{r_j \in {\cal N}(s,t)} \left( 1- \sum_{k=1}^m a_k(r_j) \right).$$
Next, let us proceed to the induction step, 
assuming that the first inequality in \eqref{G1} holds for $t \in [s+n-1,s+n)$.  
For $t\in [s+n,s+n+1)$,  
\begin{align*}
x(t) & \leq x(s+n) =  x(s+n-1) - \sum_{k=1}^m a_k(s+n-1) x(h_k(s+n-1)) \\ & \leq x(s+n-1) - \sum_{k=1}^m a_k(s+n-1) x(s+n-1)
\\
& = x(s+n-1) \left( 1- \sum_{k=1}^m a_k(s+n-1) \right) 
\\ & \leq x(s) \prod_{r_j \in {\cal N}(s,s+n-1)} \left( 1- \sum_{k=1}^m a_k(r_j) \right)
\left( 1- \sum_{k=1}^m a_k(s+n-1) \right) \\
& = x(s) \prod_{r_j \in {\cal N}(s,s+n)} \left( 1- \sum_{k=1}^m a_k(r_j) \right)
= x(s) \prod_{r_j \in {\cal N}(s,t)} \left( 1- \sum_{k=1}^m a_k(r_j) \right),
\end{align*}
by the definition of $\cal N$ in \eqref{sets}.

The proof of the second inequality in \eqref{G1} is similar, starting from $t-\ell \in [s,s+1)$ for $\ell \in {\mathbb N}$ and proceeding further to $t-\ell, t-\ell+1,
\dots, t$ by induction.
\qed

Estimate \eqref{G1} allows to prove the main result of this section.

\begin{guess}
\label{th_nonincrease}
Let
\begin{equation}
\label{oscil_1}
\limsup_{t \to \infty} \sum_{u_{\ell} \in {\cal N}(g(t),t)} \sum_{k=1}^m a_k\left( u_{\ell} \right) \prod_{r_j \in {\cal N} (h_k(u_{\ell}), 
g(t)
)}
\left( 1- \sum_{k=1}^m a_k(r_j) \right) > 1,
\end{equation} 
where $g$ is defined in \eqref{phi_k}.
Then \eqref{6} has no positive non-increasing (negative non-decreasing) solutions.
\end{guess}
{\bf Proof.}
Let $x(t)$ be a positive non-increasing solution. By \eqref{oscil_1}, there is an infinite number of $t$ such that the expression
under $\limsup$ in \eqref{oscil_1} exceeds one, and $t$ can be chosen arbitrarily large. Let us assume that $t$ is such a number,
and all $h_k(h_j(s))>0$ for any $s>t$. Assume that $d \in {\mathbb N}_0$ is the largest number such that $g(t)+ d \leq  t$.
By the choice of $t$
\begin{equation}
\label{auxil1}
\sum_{\ell=0}^{d-1} \sum_{k=1}^m a_k\left( g(t)+\ell \right) \prod_{r_j \in {\cal N} (h_k(g(t)+\ell), g(t))}
\left( 1- \sum_{k=1}^m a_k(r_j) \right) > 1.
\end{equation}

Rewriting \eqref{6} at all the points $t=g(t)$, $g(t)+1$, $\dots $, $t= g(t)+d-1$,
 $t \in {\cal N}(g(t),t)$, we get
\begin{equation}
\label{stream}
\begin{array}{l}
\displaystyle
x(g(t))=x(g(t)+1)+ \sum_{k=1}^m a_k(g(t)) x(h_k(g(t))), 
\\
\displaystyle
x(g(t)+1)=x(g(t)+2)+ \sum_{k=1}^m a_k(g(t)+1) x(h_k(g(t)+1)), \dots,
\\
\displaystyle
x(g(t)+d-1)=x(g(t)+d)+ \sum_{k=1}^m a_k(g(t)+d-1) x(h_k(g(t)+d-1))
\end{array}
\end{equation}
By Lemma~\ref{lem_flow}, 
\begin{equation}
\label{stream1}
x\left( h_k(g(t))+\ell \right) \geq x(g(t)) \prod_{r_j \in {\cal N}(h_k(g(t)+\ell),g(t))} \left( 1- \sum_{k=1}^m a_k(r_j) \right),
~~\ell=0, \dots, d-1.
\end{equation}
Summing up the $d$ equalities in \eqref{stream} and substituting inequalities \eqref{stream1}, we obtain, as 
$g(t)+d  \leq t$ and $x(g(t)+d) \geq x(t)$,
\begin{align*}
x(g(t)) - x(t) & \geq   x(g(t)) - x(g(t)+d) 
= \sum_{\ell=0}^{d-1} \sum_{k=1}^m a_k\left( g(t)+\ell \right) x \left( h_k(g(t)+\ell) \right)
\\ & \geq x(g(t)) 
\sum_{\ell=0}^{d-1} \sum_{k=1}^m a_k\left( g(t)+\ell \right) \prod_{r_j \in {\cal N} (h_k(g(t)+\ell), g(t))}
\left( 1- \sum_{k=1}^m a_k(r_j) \right) > x(g(t)),
\end{align*}
which by \eqref{auxil1} immediately implies $x(t)<0$. This contradicts to the assumption that $x(t)>0$, 
hence there is no positive non-increasing solution.
The case of negative non-decreasing solutions is considered similarly.
\qed

\begin{remark}
In \eqref{oscil_1}, we can use ${\cal M}$ instead of ${\cal N}$, getting another oscillation condition.
\end{remark}

\section{Discussion}
\label{sec5}

Let us note that existence of a positive solution is not sufficient for existence of a positive non-increasing solution.

\begin{example}
\label{example4}
The equation
\begin{equation}
\label{ex4eq1}
x(t+1)-x(t)+\left( \frac{1}{2}- \frac{1}{2} \{ t \} \right)x(t) = 0
\end{equation}
can be rewritten as 
$x(t+1)= \left( \frac{1}{2} + \frac{1}{2} \{ t \} \right)x(t)$.
Obviously 
$$
x(t)=1, ~~t \in [0,1], \quad x(t)= \left( \frac{1}{2} + \frac{1}{2} \{ t \} \right)^{\lfloor t \rfloor}
=\left( \frac{1}{2} + \frac{1}{2} \{ t \} \right)^n, ~~t \in [n,n+1).
$$
is a positive piecewise continuous solution of \eqref{ex4eq1}.
Note that it is a simple case of an equation with continuous time and integer delay,
existence of positive solutions for such equations was studied in detail in the recent paper \cite{C_Gyori_PS}
in terms of the generalized characteristic equation, which in the case of \eqref{ex4eq1}
has the solution $\displaystyle \lambda(t) =  \frac{1}{2} + \frac{1}{2} \{ t \}$.
\vspace{2mm}

On the other hand, equation \eqref{ex4eq1} 
has no positive non-increasing solutions.
Let us assume the contrary that $x(t)$ is a positive non-increasing solution, and denote the ratio by $\alpha:=x(0)/x(0.5)$,
$\alpha\in [1,\infty)$. Take $\displaystyle n>\frac{\ln \alpha}{\ln 1.5}$,
then
$\displaystyle
\frac{x(n)}{x(n+\frac{1}{2})} = \frac{\alpha 0.5^n}{0.75^n}= \frac{\alpha}{1.5^n}<1,
$
and $x(n+\frac{1}{2})>x(n)$ for $\displaystyle n>\frac{\ln \alpha}{\ln 1.5}$ large enough, which contradicts to the assumption that $x(t)$ is non-increasing.
\end{example}

The results of the present paper give rise to the following questions.

\begin{enumerate}
\item
According to Example~\ref{example4}, existence of a positive solution of the generalized characteristic equation
\cite{C_Gyori_PS} guarantees existence of a positive, but not necessarily positive non-increasing solutions.
Is it possible to develop sufficient conditions, in similar terms, for existence of positive non-increasing solutions of \eqref{6}?

\item
Can the results of \cite{Zhang1999}, where a linearized oscillation theory is developed for constant delays and coefficients, be extended to the equation 
\begin{equation*}
x(t+1)-x(t)+ \sum_{k=1}^m a_k(t)f_k(x(h_k(t)))=0, ~~a_k(t) \geq 0, ~~h_k(t) \leq t, ~~ k=1, \dots, m,~ t \geq 0,
\end{equation*}
where  $u f_k(u)>0$, $u \neq 0$ and $\displaystyle \lim_{u\to 0} \frac{f_k(u)}{u}=1$, $k=1, \dots, m$?
\item
Will sharper than in \eqref{G1}  estimates of the rate of decay of $x$, using the 
iterative procedure similar to \cite{BCS1,BCS2}, lead to a substantial improvement of  \eqref{oscil_1}?
\end{enumerate}

\section*{Acknowledgment}

Both authors were partially supported by the NSERC research grant RGPIN-2015-05976.
The authors are grateful to anonymous reviewers whose valuable comments contributed to the presentation of the paper.

\end{document}